\input amssym.def
\input amssym
\magnification=1200
\parindent0pt
\hsize=16 true cm
\baselineskip=13  pt plus .2pt
$ $
\def\M{{\cal M}}
\def\T{{\cal T}}
\def\F{{\cal F}}
\def\Z{\Bbb Z}

\def\A{\Bbb A}
\def\S{\Bbb S}

\centerline {\bf  A note on minimal finite quotients of mapping class
groups}

\bigskip \bigskip

\centerline {Bruno P. Zimmermann}

\bigskip

\centerline {Universit\`a degli Studi di Trieste}
\centerline {Dipartimento di
Matematica e Informatica}
\centerline {34100 Trieste, Italy}
\centerline
{zimmer@units.it}

\vskip 1cm

Abstract. {\sl  We prove that the minimal nontrivial finite quotient group of
the mapping class group $\M_g$ of a closed orientable surface of genus $g$ is
the symplectic group ${\rm PSp}_{2g}(\Z_2)$, for $g=3$ and 4 (this might
remain true, however, for arbitrary genus $g>2$). We discuss also some results
for arbitrary genus $g$. }

\bigskip \bigskip

{\bf 1. Introduction}

\medskip

It is an interesting but in general difficult problem to classify the
finite quotients (factor groups) of certain geometrically significant infinite
groups. This becomes particularly attractive if the group in question is perfect
(has trivial abelianization) since in this case  each finite quotient projects
onto a minimal quotient which is a nonabelian finite simple
group, and there is the well-known list of the finite simple groups (always
understood to be nonabelian in the following).

\medskip

As an example, the finite quotients of the Fuchsian triangle group of type
(2,3,7) (two generators of orders two and three whose product has order
seven) are the so-called Hurwitz groups, the groups of orientation-preserving
diffeomorphisms of maximal possible order $84(g-1)$ of a closed orientable
surface of genus $g$. There is a rich literature on the classification of the
Hurwitz groups, and in particular on the most significant case of
simple Hurwitz groups; the smallest Hurwitz group is the projective linear or
linear fractional group ${\rm PSL}_2(7)$ of order 168, acting on Klein's
quartic of genus three.

\medskip

One of the most interesting groups in topology is the mapping class group
$\M_g$ of a closed orientable surface $\F_g$ of genus $g$ which is the group of
orientation-preserving homeomorphisms of $\F_g$ modulo the subgroup of
homeomorphisms isotopic to the identity; alternatively, it is the
"orientation-preserving" subgroup of index two of the outer automorphism group
Out($\pi_1(\F)$) of the fundamental group. It is well-known that $\M_g$ is a
perfect group, for $g \ge 3$ ([Po]). By abelianizing the fundamental group
$\pi_1(\F)$ and reducing coefficients modulo a positive integer $k$, we get
canonical projections
$$\M_g  \to  {\rm Sp}_{2g}(\Bbb Z)  \to {\rm Sp}_{2g}(\Z_k) \to  {\rm
PSp}_{2g}(\Z_k)$$

of the mapping class group $\M_g$  onto the symplectic group
${\rm Sp}_{2g}(\Bbb Z)$ and the finite projective symplectic groups
${\rm PSp}_{2g}(\Z_k)$ (see [N]); we note that, for primes $p$ and $g \ge 2$,
${\rm PSp}_{2g}(\Z_p)$ is a simple  group with the only exception
of ${\rm PSp}_4(\Z_2)$ which is isomorphic to the symmetric group ${\rm S}_6$.
The kernel of the surjection $\M_g  \to  {\rm Sp}_{2g}(\Bbb Z)$ is the Torelli
group $\T_g$ of all mapping classes which act trivially on the first homology
of the surface $\F_g$.

\medskip

It is well-known that the symplectic groups ${\rm Sp}_{2g}(\Bbb Z)$ and the
linear groups ${\rm SL}_n(\Bbb Z)$ are perfect, for $g \ge 3$ resp. $n \ge 3$.
As a consequence of the congruence subgroup property for these groups, the
following holds ($p$ denotes a prime number):

\bigskip

{\bf Theorem 1.} {\sl

i) For $n \ge 3$, the finite simple quotients of the linear group ${\rm
SL}_n(\Bbb Z)$ are the linear groups ${\rm PSL}_n(\Z_p)$.

ii) For $g \ge 3$, the finite simple quotients of the symplectic group ${\rm
Sp}_{2n}(\Bbb Z)$ are the symplectic groups ${\rm PSp}_{2n}(\Z_p)$.}

\bigskip

Theorem 1 will be proved in section 4.
For the case of mapping class groups, the following is the main result of
the present note.

\bigskip

{\bf Theorem 2.} {\sl  For $g=3$ and 4, the minimal nontrivial finite quotient
group of the mapping class group $\M_g$ of genus $g$ is the
symplectic group ${\rm PSp}_{2g}(\Z_2)$.}

\bigskip

We note that the order of ${\rm PSp}_{2g}(\Z_2) = {\rm Sp}_{2g}(\Z_2)$ is
$2^{g^2}(2^2-1)(2^4-1) \ldots (2^{2g}-1)$, so for $g=3$ and 4 the orders are
1.451.520 and 47.377.612.800; these orders grow very fast, in fact
exponentially with $g^2$, whereas the orders of the finite subgroups of $\M_g$
grow only linearly with $g$ (bounded above by $84(g-1)$, see [Z, Theorem
2.1]).

\bigskip

Theorem 2 raises more questions than it answers, e.g. (even for $g=3$ this
seems to be unknown):

\medskip

- which are the finite simple quotients of $\M_g$?

- what is the minimal index of any subgroup of $\M_g$?

\medskip

Nevertheless, the proof of the Theorem appears nontrivial and interesting:
considering for $g=3$ and 4 the list of the finite simple groups of order
less than that of  ${\rm PSp}_{2g}(\Z_2)$ we are able to exclude all of
them by considering certain finite subgroups of $\M_g$ which must inject. Since
such problems are of a strongly computational character, some lists of simple
groups and case-by-case analysis seem unavoidable; also, since there does not
seem to be much relation between mapping class groups for different genera $g$,
it may be difficult to generalize Theorem 2 for arbitray $g$ (if it remains
true).  Concerning the case of  genus three, the only simple groups, different
from ${\rm PSp}_6(\Z_2)$ and of an order smaller than the order 4.585.351.680
of
${\rm PSp}_6(\Z_3)$, which we cannot exclude at moment as a quotient of $\M_3$
are the groups $^3{\rm D}_4(2)$, ${\rm M}^{\rm c}{\rm L}$ and  ${\rm
PSU}_3(\Z_{17})= {\rm U}_3(17)$ (in the notation of [C]).  For the construction
of finite quotient groups of mapping class groups, see also [Sp] and
[T]: most of these groups are again closely related to the symplectic groups
${\rm PSp}_{2g}(\Z_k)$.

\medskip

In section 3, we prove Theorem 2 for the easier case $g=3$. In section 4, we
discuss some results for arbitrary genus $g$ and then deduce Theorem 2 for the
case $g=4$; we prove also the following Theorem (we note that, by a result of
Wiman, for $g \ge 2$ the maximal order of a cyclic subgroup of $\M_g$ is
$4g+2$).

\bigskip

{\bf Theorem 3.} {\sl For $g \ge 3$, let $\phi:\M_g \to G$ be a surjection of
$\M_g$ onto a  finite simple group $G$. Then $G$ is isomorphic to a
symplectic group ${\rm PSp}_{2g}(\Z_p)$, or $G$ has an element
of order $4g+2$.}

\bigskip

{\bf 2. Proof of  Theorem 2 for $g=3$}

\bigskip

Let $G$ be a finite group of orientation-preserving diffeomorphisms of a closed
surface $\F_g$ of genus $g>1$. Then the quotient $\F_g/G$ is a closed
2-orbifold: the underlying topological space is again a closed surface of some
genus $\bar g$, and there are finitely many branch points of orders $n_1,
\ldots, n_k$; we will say that the $G$-action is of type  $(\bar g; n_1,
\ldots, n_k)$.

\bigskip

One can give the surface $\F_g$ a hyperbolic or complex structure such
that $G$ acts by isometries resp. by conformal maps of the Riemann surface, by
just uniformizing the quotient orbifold $\F_g/G$ by a Fuchsian group of
signature  $(\bar g; n_1, \ldots, n_k)$ (see e.g. [ZVC]). Then this Fuchsian
group is obtained as the group of all lifts of elements of $G$ to the universal
covering of $\F_g$ (which is the hyperbolic plane), and there is a surjection of
this Fuchsian group onto $G$ whose kernel is the universal covering group of
the surface. We will say in the following that the finite $G$-action is given by
a  surjection of a Fuchsian group of type or signature $(\bar g; n_1, \ldots,
n_k)$ onto  $G$.

\bigskip

By [FK, Theorem V.3.3], every conformal map of a closed Riemann surface of
genus
$g>1$ which induces the identity on the first homology is the identity. In
particular, every finite group of orientation-preserving diffeomorphisms of a
closed surface of genus $g>1$ injects into the mapping class group $\M_g$ and
its quotient, the symplectic group  ${\rm PSp}_{2n}(\Bbb Z)$, and we will speak
in the following of the finite group of mapping classes $G$ of $\F_g$, of type
$(\bar g; n_1, \ldots, n_k)$, determined by a surjection
$$(\bar g; n_1, \ldots,n_k) \to G$$
 of a Fuchsian group of type  $(\bar g; n_1,\ldots, n_k)$ onto $G$.

\bigskip

As an example, the Hurwitz action of the linear fractional group
${\rm PSL}_2(7)$ on the surface $\F_3$ of genus three (or Klein's quartic) is
determined by a surjection (unique up to conjugation in
${\rm PGL}_2(7)$)
$$(2,3,7)  \to  {\rm PSL}_2(7)$$
of the triangle group (0;2,3,7)=(2,3,7) onto the linear fractional group
${\rm PSL}_2(7)$, so  this defines a subgroup
${\rm PSL}_2(7)$ of the mapping class group $\M_3$.

\bigskip

We will consider in the following some finite subgroups of $\M_3$, represented
by finite groups of diffeomorphisms of a surface of genus three, or
equivalently by surjections from Fuchsian groups. For a convenient list and a
classification of the finite groups acting on a surface of genus three, see
[Br].

\bigskip

Up to conjugation, $\F_3$ has three orientation-preserving involutions which
are of types (1;2,2,2,2) = (1;$2^4$), (0;$2^8$) (a "hyperelliptic
involution") and (2;-) (a free involution). For general genus $g$, the following
is proved in [MP].

\bigskip

{\bf Proposition 1.} ([MP])

{\sl i)  If $g \ge 3$ is odd,  any
involution of type (${g-1\over2};2,2,2,2)$ normally generates $\M_g$.

ii) If $g \ge 4$ is even, any involution of type ($g\over2$;2,2)
normally generates $\M_g$.}

\bigskip

In particular,  the
cyclic group $\Z_2$ of order two of $\M_3$ generated by an involution
of type (1;$2^4$) normally generates $\M_3$ and hence maps nontrivially under
any nontrivial  homomorphism $\phi: \M_3  \to G$. On the other hand, we note
that the mapping class represented by an involution of type ($2^8$) = (0,$2^8$)
lies in kernel of the canonical surjection
$$\M_3 \to {\rm PSp}_{6}(\Bbb Z) \to {\rm PSp}_{6}(\Z_2).$$

We consider now the Hurwitz action of  ${\rm PSL}_2(7)$ on $\F_3$ defined by a
surjection  $\pi:(2,3,7)  \to {\rm PSL}_2(7)$ and realizing ${\rm PSL}_2(7)$
as a subgroup of $\M_3$.  Up to conjugation, ${\rm PSL}_2(7)$ contains a unique
subgroup $\Z_2$, and the preimage $\pi^{-1}(\Z_2)$ in the
triangle group (2,3,7) is a Fuchsian group of signature (1;$2^4$)  (since the
subgroup $\Z_2$ has index four in its normalizer in ${\rm PSL}_2(7)$ which is a
dihedral group of order eight). Since ${\rm PSL}_2(7)$ is a simple group and an
involution of type (1;$2^4$) normally generates $\M_3$, we have:

\bigskip

{\bf Lemma 1.} {\sl Every nontrivial group homomorphism $\phi: \M_3 \to G$
injects ${\rm PSL}_2(7)$.}

\bigskip

Remark.  The preimage $\pi^{-1}(\Z_3)$ of the subgroup of order three of ${\rm
PSL}_2(7)$ (unique up to conjugation) is a Fuchsian group of type (1;3,3), and
$\pi^{-1}(\Z_7)$ is a triangle group of type (7,7,7) (the normalizer of $\Z_3$
is dihedral of order 6, that of $\Z_7$ is the subgroup of ${\rm PSL}_2(7)$
represented  by all upper triangular matrices which has order 21). Since  ${\rm
PSL}_2(7)$ is simple, the corresponding subgroups $\Z_3$ and $\Z_7$ of $\M_3$
inject under any nontrivial $\phi$.

\bigskip

Now, for the proof of Theorem 2, suppose that $\phi: \M_3 \to G$ is a
surjection onto a nontrivial finite group $G$ of order less than that
of  ${\rm PSp}_{6}(\Z_2)$; since $\M_3$ is perfect, we can assume that $G$ is a
finite nonabelian simple group. By Lemma 1, $G$ has  a subgroup ${\rm
PSL}_2(7)$. The nonabelian simple groups smaller than
${\rm PSp}_{6}(\Z_2)$ and having a subgroup ${\rm PSL}_2(7)$ are the following
(in the notation of [C, p. 239ff] to which we refer for the simple groups of
small order as well as their subgroups):

$${\rm L}_2(7), \;\;  \A_7,  \;\;  {\rm U}_3(3),  \;\; \A_8,  \;\;
{\rm L}_3(4), \;\;  {\rm L}_2(49), \;\; {\rm U}_3(5),  \;\;  \A_9,  \;\;
{\rm M}_{22},  \;\; {\rm J}_2,$$

where ${\rm L}_n(p^r) = {\rm PSL}_n({\rm GF}(p^r))$ denotes a linear group over
the Galois field with $p^r$ elements,
${\rm U}_n(p) = {\rm PSU}_n(\Z_p) = {\rm PSU}_n({\rm GF}(p))$ a unitary and
$\A_n$ an alternating group; ${\rm M}(22)$ is a Mathieu group and ${\rm J}_2$
the second Janko or Hall-Janko group.

\medskip

Now each of these groups does not have simultaneously elements of order
8, 9, or 12 (see [C]), hence the proof of  Theorem 2 follows from the
following:

\bigskip

{\bf Lemma 2.} {\sl There are cyclic subgroups of orders 8, 9 and 12  of
$\M_3$ which every nontrivial homomorphism $\phi: \M_3 \to G$ injects.}

\medskip

{\it Proof.} i) A subgroup $\Z_8$ of $\M_3$ is defined by
a surjection $\pi:(4,8,8) \to \Z_8$. The preimage $\pi^{-1}(\Z_2)$ is a Fuchsian
group of type (1;$2^4$) which hence defines a subgroup $\Z_2$ of $\M_3$ which
normally generates $\M_3$.  (See also [Sn] for the determination of the
signature of a subgroup of a Fuchsian group.)

\medskip

ii) A subgroup $\Z_9$ of $\M_3$ is defined by a surjection $\pi:(3,9,9) \to
\Z_9$, and the preimage $\pi^{-1}(\Z_3)$ gives a subgroup $\Z_3$ of $\M_3$ of
type ($3^5$) (we note that, up to conjugation, there are exactly two periodic
diffeomorphisms of order three of $\F_3$, of types  ($3^5$) and
(1;3,3)). Suppose, by contradiction, that $\phi$ is trivial on the subgroup
$\Z_3$.

\medskip

We consider a subgroup ${\rm SL}_2(3)$ of $\M_3$ defined by a surjection
$\pi:(3,3,6)) \to {\rm SL}_2(3)$; the linear group
${\rm SL}_2(3)$ of order 24 is isomorphic to the binary tetrahedral group
$\A_4^*$ and is a semidirect product $Q_8 \ltimes \Z_3$. The preimage
$\pi^{-1}(\Z_3)$ defines a subgroup $\Z_3$ of $\M_3$ of type ($3^5$) which, by
hypothesis, is mapped trivially by  $\phi$.  Now it follows easily that $\phi$
has to be trivial on the whole subgroup $Q_8 \ltimes \Z_3$, and in particular
on the unique cyclic subgroup $\Z_2$ of order two of $Q_8$ which is
of type (1;$2^4$). Since $\Z_2$ normally generates $\M_3$, the homomophism
$\phi$ is trivial.

\medskip

iii) A subgroup $\Z_{12}$ of $\M_3$ is defined by a
surjection $\pi:(3,4,12) \to \Z_{12}$; now $\pi^{-1}(\Z_2)$ is of type
(1;$2^4$), and  $\pi^{-1}(\Z_3)$ of type ($3^5$). Since $\phi$ is nontrivial,
it cannot be trivial on the subgroup $\Z_2$ of $\M_3$ of type (1;$2^4$). On the
other hand, if $\phi$ is trivial on the subgroup $\Z_3$ of type ($3^5$) then one
concludes as in ii) that $\phi$ is trivial.

\bigskip

This concludes the proof of Lemma 1 and also of the case $g=3$ of Theorem 2.

\bigskip

Remark. The unitary group ${\rm U}_3(3)$
can be excluded also by considering a quaternion subgroup $Q_8$ of order
eight of $\M_3$ defined by a surjection $\pi:(1;2) \to Q_8$. Again, the preimage
$\pi^{-1}(\Z_2)$ of the unique subgroup $\Z_2$ of $Q_8$ has signature
(1;$2^4$) and defines a subgroup $\Z_2$ of $\M_3$. However ${\rm U}_3(3)$ (whose
Sylow 2-subgroup is a wreathed product $(\Z_4 \times \Z_4) \ltimes \Z_2$) has no
subgroup $Q_8$, so $\phi$ maps $\Z_2$ trivially and hence all of $\M_3$.

\bigskip

{\bf 3. Some results for arbitrary genus; proof of Theorem 2 for $g=4$}

\bigskip

The following is the main result of [Pa].

\bigskip

{\bf Proposition 2.} ([Pa])

{\sl For $g \ge 3$, the index of any proper subgroup of $\M_g$ is larger than
$4g+4$; equivalently, there are no surjections of $\M_g$ onto an alternating
group  $\A_n$, or onto any transitive subgroup of $\A_n$, if  $3 \le n \le
4g+4$.}

\bigskip

It would be interesting to know if there exists any surjection $\phi:\M_G \to
\A_n$, for $n > 4g+4$.

\bigskip

The maximal order of a cyclic subgroup of $\M_g$ is $4g+2$, for any $g>1$, and
such a maximal subgroup $\Z_{4g+2}$ is generated by a diffeomorphism of type
(2,$2g+1,4g+2$); the subgroup $\Z_2$ of $\Z_{4g+2}$ is generated by a
hyperelliptic involution of type (0;$2^{2g+2}$). The following result is proved
in [HK].

\bigskip

{\bf Proposition 3.} ([HK]) {\sl Let $g \ge 3$.

i) Let $h$ be an orientation-preserving diffeomorphism
of maximal order $4g+2$ of $\F_g$. If $1 \le k \le 2g$ then
$h^k$ normally generates $\M_g$.

ii) The normal subgroup of $\M_g$ generated by the
hyperelliptic involution $h^{2g+1}$ contains the Torelli group $\T_g$ as a
subgroup of index two and is equal to the kernel of the canonical projection
$\M_g  \to  {\rm Sp}_{2g}(\Z)  \to {\rm PSp}_{2g}(\Z) = {\rm Sp}_{2g}(\Z)/\{\pm
I\}$.

iii) Let $G$ be a group without an element of order $g-1$, $g$ or $2g+1$. Then
any homomorphism $\phi:\M_g \to G$ is trivial.}

\bigskip

Note that i) and ii) of Proposition 3 combined with Theorem 1 imply Theorem 3.

\medskip

We consider the case $g=4$ now.

\bigskip

{\bf Lemma 3.} {\sl Let $\phi:\M_4 \to G$ be a surjection onto a finite
simple group $G$.

i) The symmetric group $\S_5$ is a subgroup of $G$.

ii) Either $G$ has elements of orders 10, 16 and 18, or $G$ is isomorphic to
a symplectic group ${\rm PSp}_8(\Z_p)$.}

\bigskip

{\it Proof.}   The mapping class group $\M_4$ has a subgroup $\A_5$ of type
(2,5,5) and a subgroup $\S_5$ of type (2,4,5). An involution in a subgroup
$\A_5$ of $\M_4$ defines a subgroup $\Z_2$ of type (2;2,2); by Proposition 1,
such a subgroup $\Z_2$ normally generates $\M_4$, and  hence $\phi$
injects $\A_5$, $\S_5$ and also their subgroups $\Z_5$ which are of type
(0;$5^4$).  Now
$\M_4$ has a subgroup $\Z_{10}$ of type (5,10,10), and since its subgroups
$\Z_5$ and $\Z_2$ are of type (0;$5^4$) and (2;2,2) and hence
inject, also  $\Z_{10}$  injects.

\medskip

Also, $\M_4$ has a subgroup $\Z_{16}$ of type (2,16,16) whose
subgroup $\Z_2$ is of hyperelliptic type (0;$2^{10}$), and a maximal
cyclic subgroup $\Z_{18}$ of type (2,9,18). Lemma 3ii) is now a
consequence of Proposition 3 and Theorem 1.

\bigskip

{\it Proof of Theorem 2 for the case $g=4$}.

\medskip

Let $\phi: \M_4 \to G$ be a surjection onto a finite simple group $G$. Suppose
that the order of $G$ is less than the
order 47.377.612.800 of ${\rm PSp}_8(\Z_2)$; see [C,p.239ff] for a list of
these groups. The alternating groups of such orders are excluded by
Proposition  2 since they have subgroups of index $\le 20$.  The linear
groups  ${\rm PSL}_2(p^r)$ in dimension two are excluded by Lemma 3 since,
with the exceptions of ${\rm PSL}_2(5^2)$ and ${\rm PSL}_2(5^4)$, they have no
subgroups $\S_5 \cong {\rm PGL}_2(5)$. All remaining groups in the list can be
excluded case by case by considering the possible orders of elements in each
of these groups (see [C] for the character tables of most of these groups;
the group theory package GAP can also be used to create the conjugacy classes
and the orders of the elements of these group). It is easy to see then
that none of these groups has simultaneously elements of orders 10, 16 and 18
(in some cases it may be helpful also to consider a Sylow 2-subgroup of $G$).
Applying Lemma 3 again completes now the proof that the smallest (simple)
quotient group of $\M_4$ is indeed the symplectic group ${\rm PSp}_8(\Z_2)$.

\bigskip

{\bf 4. Proof of  Theorem 1}

\medskip

i) Let $\phi: {\rm SL}_n(\Bbb Z) \to G$  be a surjection onto a finite simple
group $G$.  By the congruence subgroup property for linear groups in
dimensions $n>2$ (which holds also for symplectic groups, see [M],[BMS]), the
kernel of
$\phi$ contains a congruence subgroup, i.e. the kernel of a canonical
projection ${\rm SL}_n(\Bbb Z) \to {\rm SL}_n(\Z_k)$, for some
positive integer $k$, and hence $\phi$ induces a surjection
$\psi: {\rm SL}_n(\Z_k) \to G$ (see [N,II.21]).

\medskip

If $k = p_1^{r_1} \ldots p_s^{r_s}$ is the prime decomposition,
$${\rm SL}_n(\Z_k) \; \cong \; {\rm SL}_n(\Z_{p_1^{r_1}}) \times \ldots \times
{\rm SL}_n(\Z_{p_s^{r_s}})$$ (see [N,Theorem VII.11]). Now the restriction of
$\psi: {\rm SL}_n(\Z_k) \to G$ to some factor ${\rm SL}_n(\Z_{p_i^{r_i}})$ has
to be nontrivial; since $G$ is simple, this gives some surjection $\psi: {\rm
SL}_n(\Z_{p^r}) \to G$.

\medskip

Let $K$ denote the kernel of the canonical surjection ${\rm SL}_n(\Z_{p^r}) \to
{\rm SL}_n(\Z_p)$, so $K$ consists of all matrices in ${\rm SL}_n(\Z_{p^r})$
which are congruent to the identity matrix $I$ when entries are taken modulo
$p$. By performing the binomial expansion of $(I + pA)^{p^{r-1}}$ one checks
easily that $K$ is a $p$-group, in particular $K$ is solvable.  Then also  the
kernel $K_0$ of the canonical surjection from ${\rm SL}_n(\Z_{p^r})$ to the
central quotient ${\rm PSL}_n(\Z_p)$ of ${\rm SL}_n(\Z_p)$ is solvable. Since
$G$ is simple, $\psi$ maps $K_0$ trivially and induces a surjection
from ${\rm PSL}_n(\Z_p)$ onto $G$; since  $n>2$, ${\rm PSL}_n(\Z_p)$ is simple
and this surjection is an isomorphism.

\medskip

ii)  By [N,Theorem VII.26],
$${\rm Sp}_{2n}(\Z_k) \; \cong \; {\rm Sp}_{2n}(\Z_{p_1^{r_1}}) \times \ldots
\times {\rm Sp}_{2n}(\Z_{p_s^{r_s}}),$$
and the proof is then analogous to the first case.

\bigskip  \bigskip

\centerline {\bf References}

\bigskip

\item {[BMS]} H.Bass, J.Milnor, J.P.Serre, {\it The congruence subgroup
property for $SL_n$ ($n \ge 3$) and  $SP_{2n}$ ($n \ge 2$).}    Inst.
Hautes Etudes Sci. Publ. Math. 33, 59-137 (1967)

\item {[Br]} S.A.Broughton, {\sl Classifying finite group actions on
surfaces of low genus.} J. Pure Appl. Algebra 69,  233-270  (1990)

\item {[C]} J.H.Conway, R.T.Curtis, S.P.Norton, R.A.Parker, R.A.Wilson,
Atlas of Finite Groups. Oxford University Press 1985

\item {[FK]} W.J.Harvey, M.Korkmaz, {\sl Homomorphisms from mapping class
groups}. Bull. London Math. Soc. 37, 275-284 (2005)

\item {[HK]} H.M.Farkas, I.Kra, {\sl Riemann surfaces}. Second edition.
Graduate Texts in Mathematics 71, Springer 1991

\item {[MP]} J.McCarthy, A.Papadopoulos {\it Involutions in surface mapping
class groups.} L'Enseig. Math. 33, 275-290  (1987)

\item {[Me]} J.Mennicke, {\it Finite factor groups of the unimodular groups.}
Ann. Math. 81, 31-37 (1965)

\item {[N]} M.Newman, {\it Integral Matrices.} Pure and Applied Mathematics
Vol.45,  Academic Press 1972

\item {[Pa]} L.Paris, {\it Small index subgroups of the mapping class groups.}
Preprint (electronic version under arXiv:math.GT/0712.2153v1)

\item {[Po]} J. Powell, {\it Two theorems on the mapping class group of a
surface.}  Proc. Amer. Math. Soc. 68, 347-350  (1978)

\item {[Sn]}  D.Singerman,  {\it  Subgroups of Fuchsian groups and finite
permutation groups.}   Bull. London Math. Soc. 2, 319-323  (1970)

\item {[Sp]}  P.L.Sipe,  {\it  Some finite quotients of the mapping class group
of a surface.} Proc. Amer. Math. Soc. 97,  515-524  (1986)

\item {[T]}  F.Taherkhani,  {\it  The Kazhdan property of the mapping class
group of closed surfaces and the first cohomology of its cofinite subgroups.}
Experimental Math. 9,  261-274  (20007

\item {[Z]}  B.Zimmermann,  {\it  Lifting finite groups of outer
automorphisms of free groups, surface groups and their abelianizations.} Rend.
Istit. Mat. Univ. Trieste 37,  273-282  (2005) (electronic version in
arXiv:math.GT/0604464)

\item {[ZVC]} H.Zieschang, E.Vogt, H.-D.Coldewey, Surfaces and planar
discontinuous groups. Lecture Notes in Mathematics 835, Springer, Berlin
1980

\bye